\documentclass[10pt,twoside]{article}
\usepackage{graphicx}
\usepackage{amsmath}
\usepackage{Latex-document}

\newcommand{\be}{\begin{eqnarray}}
\newcommand{\en}{\end{eqnarray}\vs 0.5 cm}
\newcommand{\non}{\nonumber}
\newcommand{\no}{\noindent}

\newcommand{\Ns}{{\bf s}}

\newcommand{\Nna}{{\bf \nabla}}

\newcommand{\NR}{{{\bf R}}}
\newcommand{\NT}{{{\bf T}}}
\newcommand{\NZ}{{{\bf Z}}}

\newcommand{\Nn}{{{\bf n}}}

\newcommand{\qq}{\begin{eqnarray}}

\newcommand{\da}{\partial}

\newcommand{\qqq}{\end{eqnarray}}

\newcommand{\tr}{\hbox{tr}}

\newcommand{\CB}{{\cal B}}

\newcommand{\hf}{{_1\over^2}}

\newcommand{\vs}{\vskip}

\def\volumeno{I}

\title{\bf Ergodicity and Mixing  for \vskip -2mm
Stochastic Partial Differential Equations\vskip 6mm}

\author{J. Bricmont\vspace*{-0.5cm}\thanks{UCL,
Physique Th\'eorique, Chemin du Cyclotron 2, B-1348, Louvain-la-Neuve,
Belgium. E-mail: bricmont@fyma.ucl.ac.be}}

\date{\vskip -8mm}

\begin{document}

\maketitle

\thispagestyle{first} \setcounter{page}{567}

\begin{abstract}

\vskip 3mm

Recently, a number of authors have investigated the conditions
under which a stochastic perturbation acting on an infinite
dimensional dynamical system, e.g. a partial differential
equation, makes the system ergodic and mixing. In particular, one
is interested in finding minimal and physically natural conditions
on the nature of the stochastic perturbation. I shall review
recent results on this question; in particular, I shall discuss
 the Navier-Stokes equation on a
two dimensional torus with a random force which is white noise
in time, and excites only
a finite number of  modes. The number of excited
modes depends on the viscosity $\nu$, and
 grows like  $\nu^{-3}$ when $\nu$
goes to zero.  This Markov process
has a  unique  invariant measure and is
 exponentially mixing in time.

\vskip 4.5mm


\end{abstract}

\vskip 12mm

\section{Introduction} \label{section 1}\setzero
\vskip-5mm \hspace{5mm }

The goal of this paper is to consider stochastic partial differential equations
and to
study  conditions on the random perturbation that imply exponential
convergence to a stationary state. In fact, one wants `minimal' conditions,
in the
following sense:  by expanding the solution in a basis of eigenfunctions of
a linear operator
associated with the PDE, one can write the latter as an infinite
dimensional system of
coupled differential equations. The question, then, is: to how many such
equations do we
need to add noise in order to make the system ergodic and mixing?

The physical motivation for this question comes from the fact that
isotropic turbulence is often mathematically
modelled by Navier Stokes equation subjected to an external
stochastic driving force which is stationary in space and time.
If the solution is expanded into Fourier modes, the driving force, which,
in the
language of physicists, acts on ``large scale", should not perturb, or
perturb very weakly,
the high modes which represent the small scale properties of the system.
So, one  would like
to show that the system becomes ergodic and mixing by adding noise to as
few modes as
possible. Obviously, this requires some detailed understanding of the
nonlinear dynamics of
the deterministic PDE.

This problem is interesting from another point of view. As we shall see
below, one can show
that all but a finite number of modes converge to equilibrium provided the
remaining
ones do. So, we can reduce ourselves to a finite dimensional problem, which
would be
standard, except for the fact that the discarded modes produce a memory
effect on the
remaining ones, so that the problem is no longer Markovian. At this point,
one introduces
techniques coming from the study  of the statistical mechanics of
one-dimensional systems
(where the unique dimension corresponds physically to space rather than
time) with ``long
range, exponentially decaying, interactions"  which have already been very
useful in the
study of SRB measures in dynamical systems (see
\cite{Si}).

At present, the best results require that the number of modes to which
noise must be added
depends on the parameters of the system, although a stronger result is
likely to hold (see
Remark 4 after Theorem 1.1).

The type of question discussed here (for the Navier-Stokes equation but also
 for other equations) has been
at the center of attention of several groups of people
(see Remark 3 after Theorem 1.1 below). In this paper, I shall try to explain,
in a simplified form, the approach followed by A. Kupiainen, R. Lefevere
and myself in \cite{bk3} (see also
\cite{bkl1, bkl2} for previous results).

\vspace*{3mm}

To be concrete, consider the stochastic Navier-Stokes equation for the
velocity field ${u}(t,x)\in \NR^2$ defined on the torus $\NT=
(\NR/2\pi \NZ)^2$:
\qq
d{ u}+(({ u}\cdot\Nna){ u}-{\nu}\Nna^2{ u}+\Nna p)dt=d{ f}
\label{ns}
\qqq
where ${ f}(t,x)$ is a Wiener process with covariance
\qq
Ef_\alpha(t,x)f_\beta(t',y)= \min\{t,t'\}
C_{\alpha\beta}(x-y)
\label{F}
\qqq
and $C_{\alpha\beta}$ is a smooth function
satisfying $\sum_\alpha\da_\alpha C_{\alpha\beta}=0$.
Equation (\ref{ns}) is supplemented with
the incompressibility condition $\Nna\cdot{ u}=0=\Nna\cdot{ f}$,
and we will also assume that
the  averages over the torus vanish:
$\int_{\NT}{ u}(0,x)=0=\int_{\NT}{ f}(t,x)$, which imply
that $\int_{\NT}{ u}(t,x)=0$
 for all times $t$.

It is convenient to change (\ref{ns}) to dimensionless variables so that $\nu$
becomes equal to one.  This is achieved by setting
$
{ u}(t,x)=\nu u'(\nu t, x).
$
Then $ u'$  satisfies
(\ref{ns}), (\ref{F}) with $\nu$ replaced
by $1$, and $C$ by
$$
C'=\nu^{-3}C.
$$
From now on, we work with such variables
and drop the primes. The dimensionless control parameter
in the problem is the (rescaled) energy injection rate
$\hf\tr\, C'(0)$ ,
customarily written as
$
({\rm Re})^3$
where ${\rm Re}$ is  the
 Reynolds number:
$$
{\rm Re}=\epsilon^{_1\over^3}\nu^{-1},
$$
and $\epsilon= \hf\tr\, C(0)$ is the energy injection rate in the
original units (for explanations of the terminology see \cite{frisch}).

In two dimensions, the incompressibility condition can be
conveniently solved by expressing the velocity field
in terms of the vorticity $\omega=\da_1u_2-\da_2u_1$. First
(\ref{ns}) implies the transport equation
\qq
d\omega+(({ u}\cdot\Nna)\omega-\Nna^2\omega)dt=db,
\label{ve}
\qqq
where $b=\da_1f_2-\da_2f_1$ has the covariance
\qq
Eb(t, x)b(t', y)=\min\{t,t'\}(2\pi)^{-1}\gamma( x- y)
\non
\qqq
with $\gamma = -2\pi\nu^{-3}\Delta\tr C$.

Next, going to the Fourier transform, $\omega_
k(t)={_1\over^{2\pi}}\int_{\NT}
e^{i k\cdot x}\omega(t, x)d x$, with $ k\in\NZ^2$; we may
express $ u$ as  $ u_k=i{_{(-k_2,k_1)}\over^{k^{2}}}
\omega_k$, and write
the vorticity equation as
\qq
d\omega(t)=F(\omega(t))dt +db(t),
\label{ee}
\qqq
where the drift is given by
\qq
F(\omega)_k=-k^2\omega_k+{_1\over^{2\pi}}\sum_{l\in\NZ^2\backslash\{{
0},k\}}
{_{k_1l_2-l_1k_2}\over^{|l|^{2}}}
\omega_{k-l}
\omega_l
\label{drift}
\qqq
and $\{b_k\}$ are Brownian
motions
 with ${\bar b_k}=b_{-k}$ and
$$
Eb_k(t)b_l(t')=\min\{t,t'\}\delta_{k ,
-l}\,\gamma_k .
$$

The dimensionless control parameter for the vorticity equation is
\qq
R=\sum_{k\in\NZ^2}
\gamma_k=2\pi\gamma(0)
\label{s3}
\qqq
which is proportional to the $\omega$ injection rate,
and also to the third power of the
 Reynolds number. One is interested
in the turbulent region, where $R$ is large; therefore,
we will always assume below, when it is convenient,
 that $R $ is sufficiently large.

For turbulence one is interested in the properties of stationary state
of
the stochastic equation (\ref{ee}) in the case of {\it smooth} forcing
(see \cite{bkl1} for some discussion of this issue) and,
ideally, one would like to consider the case
where one excites only a finite number of modes,
$$
\gamma_k\neq 0 \;,\; k^2\leq N,
$$
 with $N$ of  order of one (for that, see Remark 4 below). In this paper we
assume
that $N$ scales as
\qq
N=\kappa R,
\label{kappa}
\qqq
with $\kappa$ a constant, taken large enough.
We set all the other $
\gamma_k=0$, although this condition can easily be
relaxed. Let us denote the minimum of the covariance
by
$$
\rho=\min\{|\gamma_{k}|\;|\;|k|^2\leq N\}.
$$

Before stating our result, we need some definitions.
Let $P$ be the orthogonal projection in $H=L^2(\NT)$
to the subspace $H_s$ of functions having zero Fourier
components for $|k|^2> N$. We will write
$$
\omega=s+l
$$
with $s=P\omega$, $l=(1- P)\omega$
(respectively, the small $k$ and large $k$ parts of $\omega$).
Denote also by $H_l$ the complementary subspace (containing
the nonzero components of $l$).
$H$ is our probability space, equipped with $\CB$, the
 Borel
$\sigma$-algebra.

The stochastic equation (\ref{ee}) gives rise to a Markov
process $\omega(t)$ and we denote by
$P^t(E|\omega)$ the transition probability of this process.

The  main result of \cite{bk3} is the

\vspace{3mm}

\no{\bf Theorem 1.1.} {\it The stochastic Navier-Stokes equation
} (\ref{ee}) {\it defines a  Markov process with state space
$(H, \CB)$ and  for all $R<\infty$, $\rho>0$ it has a unique
invariant measure $\mu$ there.  Moreover,
$\forall \omega\in H$, for all Borel sets $ E\in H_s$ and
for all bounded H\"older continuous functions $F$ on $H_l$, we have,
\qq
|\int P^t(d\omega'|\omega)1_E(s')F(l')-\int\mu(d\omega')1_E(s')F(l'))
|\leq C ||F||_\alpha e^{-mt}
\label{con}
\qqq
where $C= C(\|\omega\|, R, \rho)
<\infty$, $m=m(R,\rho, \alpha)>0$, and $||F||_\alpha$ is the  H\"older
norm of exponent $\alpha$.}

\vskip 2mm

\no{\bf Remark 1.}  In \cite{bk3}, we stated, for convenience,
Theorem 1.1 by saying
that the constant $C$ in (\ref{con}) was a function of $\omega$ which was
almost
surely finite. Since  this was stressed e.g. in \cite{K2}, it is worth
remarking that
$C$ is simply a function of $\|\omega\|$ (depending also on the
parameters
$R$ and $\rho$), which is finite $\forall \omega\in H$.
 To check this, we  refer the reader to equations (86)
and (97) in \cite{bk3}.
The main reason why this bound
holds, however, lies in the fact that the only dependence of our estimates
on $\omega$
appears in Lemma 4.1 below and occurs through $\|\omega\|$.

\vskip 2mm

\no{\bf Remark 2.} In  \cite{bkl1} it was proven
 that, with probability 1, the functions on the support of the
measure constructed here are real analytic. In particular
all correlation functions of the form
$$
\int\mu(d\omega)\prod_i\nabla^{n_i}u(x_i)
$$
exist.
For further results on analyticity, see \cite{Sh, MAT2}.

\vskip 2mm

\no{\bf Remark 3.} While the existence of the invariant measure
follows with soft methods \cite{vf}, its uniqueness and the
ergodic and mixing properties of the process has been harder to
establish. With a nonsmooth forcing (meaning that the strength of
the noise, $\gamma_k$, decays only polynomially with $|k|$) this
was established in \cite{fm} and for large viscosity in
\cite{mat}. However, those results did not cover the most
physically interesting situations. The first result for a smooth
forcing was by Kuksin and Shirikyan \cite{ks} who considered a
periodically kicked system with bounded kicks (for results on
exponential
 convergence in that model, see \cite{KS1, KS3, K1, MY}).
In particular they could deal with the case where only a finite number of modes
 are excited by the noise
(the number of modes depends both on the
viscosity and the size of the kicks). In \cite{bkl2},
we proved uniqueness and exponential mixing for such
a kicked system where the kicks have a  Gaussian distribution, but we required
that there be a nonzero noise for each mode.
 An
essential ingredient in  analysis of  \cite{ks}, which was used
in \cite{bk3} and by other authors, is the Lyapunov-Schmidt
type reduction that allows
to transform the original Markov process with infinite dimensional
state space to  a non-Markovian process with
finite dimensional state space.
While the analysis of \cite{ks} was limited to bounded noise acting at
discrete times, it was
extended in \cite{KS2, KS4, K2, KS6} to cover unbounded noise
and continuous time, as well as to obtain results on the strong law of
large number and
the central limit theorem.
 The first
results on ergodicity of the system with unbounded noise and finitely many
excited modes were obtained in
\cite{EMS1, bk3} (see also \cite{EL} for applications to other equations)
and, for exponential convergence, in
\cite{bk3}, which was  also proved in \cite{MAT1}. For results on related
problems, see \cite{EH, HAI, HAI2,
MY}.

\vskip 2mm

\no{\bf Remark 4.} What one would like to obtain is a result similar to
Theorem 1.1, but
with $N$ finite, independently of $R$. An interesting result in that
direction was obtained
by Weinan E and Mattingly \cite{EMS2} who showed that, if one adds noise to
only $2$ (suitably chosen)
 modes, ergodicity holds, provided one truncates the system (\ref{ee},
\ref{drift}), by keeping only a finite,
but arbitrarily large, number of modes. This of course suggests that the $2$
stochatically perturbed modes produce an ``effective
noise" on any finite number of modes, in particular on all those with $k^2
\leq \kappa R$; then, one could
hope to combine this with the results  in \cite{EMS1, bk3} to obtain
ergodicity and mixing for the full
system. This, however, has not been done.

\vskip 2mm

\no{\bf Remark 5.} The parameters in our problem are
$R$ and $\rho$. All constants that do not depend on them
will be generically denoted by $C$ or $c$.
These constants
can vary from place to place.

\vskip 3mm

Let me now explain
 the  connection with ideas coming from
statistical mechanics.

First,  observe that,
if one neglects the nonlinear term in
(\ref{ee}-\ref{drift}), one expects $\|\omega\|$ to be of
order $R^{\hf}$, for typical realizations of the noise
($R^{\hf}$ is the typical size of the noise, and the
$-k^2
\omega_k$ term will dominate in eq. (\ref{ee})
for larger values of $\|\omega\|$). It turns out that
similar probabilistic estimates hold for the full
equation (\ref{ee}) as shown in Section 4. Now, if
$\|\omega\|$ is of size $R^{\hf}$, the $-k^2
\omega_k$ term will dominate the nonlinear
 term (which is roughly of size $\|\omega\|^2$) in eq.
(\ref{ee}), for $|k| \geq \kappa R^{\hf}$, and one can
expect that those modes (corresponding to $l$ above) will
behave somewhat like the solution of the heat equation
and, in particular, that they will converge to a
stationary state.

Thus, the first step is to express the $l$-modes in terms
of the $s$-modes at previous times. This is done in Section 2 and
produces a process for the $s$-modes that is no longer
Markovian but has an infinite memory. In statistical
mechanics, this would correspond to a system of unbounded
spins (the $s$-modes) with infinite range interactions,
with the added complications that, here, the measure is
not given in a Gibbsian form, but only through a Girsanov
formula, i.e. (\ref{s1}) below, and that time is
continuous. Hence, we  have to solve several problems: the
possibility that $\omega$ be atypically large, the long
range ``interactions", and finally, showing that a
version of the $s$-process with a suitable cutoff is
ergodic and mixing.

In Section 3, I introduce  a ``toy model", namely a process with
infinite memory, but with bounded variables, so that the problems
caused by the unprobably large values of $\|\omega\|$ does not
occur. In that model, I explain how the statistical mechanical
techniques, developed to study systems on one dimensional
lattices, can be adapted to our setting.

The large $\omega$ problem is treated in Section 4, using
probabilistic estimates developed in \cite{bkl1}, which, in
statistical mechanics, would be called stability estimates. In
Section 5, I sketch how the remaining problems are handled:
showing that the techniques explained in Section 3 can be applied
here. However, this is where several technical complications
enter, for the treatment of which I refer to \cite{bk3}. The
problem is that, even though for typical noise, hence for typical
$\omega$'s, the $l$-modes depend exponentially weakly on their
past (see Section 2), thus producing, typically,  ``interactions"
that decay exponentially fast, they may depend sensitively on
their past when the noise is large. In the language of statistical
mechanics, atypically large noise produces long range
correlations, and that is the source of many technical
difficulties. My goal here is to present the main conceptual tools
used in \cite{bk3}, putting aside those difficulties.

\section{Finite dimensional reduction} \label{section 2}
\setzero\vskip-5mm \hspace{5mm }

Using an idea of \cite{ks}, one can reduce the problem of the
study of a Markov process with infinite dimensional
state space to that of a non-Markovian process with
finite dimensional state space.

For this purpose, write the equation (\ref{ee}) for
the small and large components of $\omega$ separately:
\qq
ds(t)&=&PF(s(t)+l(t))dt +db(t)\label{seq}\\
{_d\over^{dt}}l(t)&=&(1-P)F(s(t)+l(t)) . \label{leq} \qqq The idea
of  \cite{ks} is to solve the $l$ equation for a given function
$s$, thereby defining $l(t)$ as a function of the entire history
of $s(t')$, $t'\leq t$. Then, the $s$ equation will have a drift
with memory. Let us fix some notation. For a time interval $I$, we
denote the restriction of $\omega$ (or $s$, $l$ respectively) to
$I$ by $\omega(I)$, and use the boldface notation ${\bf s}(I)$, to
contrast it with $s(t)$, the value of $s$ at a single time.
 $\|\cdot\|$ will denote the
$L^2$ norm. In \cite{bkl1} it was proven that, for any $\tau<\infty$,
there exists a set $\CB_\tau$
of Brownian paths $b\in C([0,\tau],H_{s})$ of full measure
such that, for $b\in \CB_\tau$,
(\ref{ee}) has a unique solution with
$||\omega(t)||<\infty$, $||\nabla\omega(t)||<\infty$
for all $t$ (actually, $\omega(t)$ is  real analytic).
In particular, the projections $s$ and $l$ of this solution
are in $C([0,\tau],H_{s(l)})$ respectively.

On the other hand, let us denote, given any ${\bf s}\in
C([0,\tau],H_s)$, the solution --- whose existence  will be
discussed below --- of (\ref{leq}), with  initial condition $l(0)$
by $l(t,{\bf s}([0,t]),l(0))$. More generally, given initial data
$l(t')$ at time $t'<\tau$ and ${\bf s}([t',\tau])$, the solution
of (\ref{leq}) is denoted, for $\sigma \leq \tau$, by $l(\sigma,
{\bf s}([t',\sigma]),l(t'))$ and the corresponding $\omega$ by
$\omega(\sigma, {\bf s}([t',\sigma]),l(t'))$. The existence and
key properties of those functions are given by:

\vs 2mm

\no{\bf Proposition 2.1.} {\it Let $l(0)\in H_l$ and
$s\in C([0,\tau], H_s) $ . Then $l(\cdot,{\bf s}([0,t]),l(0))\in
C([0,\tau], H_l)
\cap L^2([0,\tau], H^1_l)$, where $H^1_l= H_l\cap H^1$, and $H^1$ is
the first Sobolev space.
Moreover, given two initial conditions $l_1 , l_2$ and $t\leq \tau$
\qq
\| l (t, {\bf s}([0,t]),l_1) - l(t,{\bf s}([0,t]),l_2) \| \leq
\exp \left[- \kappa Rt
+ a\int^t_0 \| \Nna \omega_1 \|^2\right]
\| l_1-l_2\|
\label{21}
\qqq
where $a=(2\pi)^{-2}\sum |k|^{-4}$ and
$\omega_1(t)=s(t)+l_1(t,{\bf s}([0,t]),l_1)$. The solution also satisfies
\qq
l(t, {\bf s}([0,t]),l(0))=l(t, {\bf s}([\tau,t]),l(\tau, {\bf
s}([0,\tau]),l(0))).
\label{22}
\qqq
}

\vskip 2mm

\no{\bf Remark.} What this Proposition shows is that the dependence of the
function
$l$ upon its initial condition $l_i$, $i=1,2$, decays exponentially in time
(i.e. like the
solution of the heat equation), provided $\omega$ is not too large, in the
sense that
$\int^t_0 \| \Nna \omega_1 \|^2\leq cRt$, for a suitable constant $c$. As
we will see in
Section 4, this event is highly probable.

\vs 4mm

Now, if $s=P\omega$ with $\omega$ as above being the solution of
(\ref{ee}) with noise $b\in \CB_\tau$ then the $l({\bf s})$
constructed in the Proposition equals $(1-P)\omega$ and the
stochastic process $s(t)$ satisfies the reduced equation \qq
ds(t)= f(t)dt+db(t) \label{sred} \qqq with \qq f(t)=
PF(\omega(t)). \label{f} \qqq where $\omega(t)$ is the function on
$C([0,t],H_s)\times H_l$ given by \qq \omega(t)=s(t)+l(t,{\bf
s}([0,t]),l(0)). \label{oomega} \qqq (\ref{sred}) has almost
surely bounded paths and we have a Girsanov representation for the
transition probability of the $\omega$-process in terms of the
$s$-variables \qq P^t(F|\omega(0)) =\int
\mu^t_{\omega(0)}(d\Ns)F(\omega(t)) \label{girs} \qqq with \qq
\mu^t_{\omega(0)}(d\Ns)=e^{\int_0^t(f(\tau),\gamma^{-1}(ds(\tau)
-\hf f(\tau)d\tau))}\nu^t_{s(0)}(d\Ns) \label{s1} \qqq where
$\nu^t_{s(0)}$ is the Wiener measure with covariance $\gamma$ on
paths ${\bf s}={\bf s}([0,t])$ with starting point $s(0)$ and
$(\cdot,\cdot)$ the $\ell^2$ scalar product.  Define the operator
$\gamma^{-1}$ in terms of its action on the Fourier coefficients:
\qq (f,\gamma^{-1}f)=\sum_{|k|^2\leq N}|f_k|^2\gamma_k^{-1}.
\label{norga} \qqq

\vs 3mm

The Girsanov representation (\ref{girs}) is convenient since
the problem of a stochastic PDE has been reduced to that of
a stochastic process with finite dimensional state space.
The drawback is that this process has infinite memory.
In the next section, I will show how to deal with this problem in a
simplified situation.

\section{A Toy Model}
\label{section 3}
\setzero\vskip-5mm \hspace{5mm }

In order to explain the main ideas in the proof, I will consider first
a `toy model' and then explain the steps needed to control the full model.

Let us consider variables $x_t \in [0,1]$, $t \in \NZ$ about which
a set of (consistent) conditional probability densities $p(x_t |
{\bf x}_{[-\infty,t-1]})$ is given, i.e. one is given
the probability densities of the variables
$x_t$, at time $t$, given a `past history'  ${\bf x}_{[-\infty,t-1]}$,
where we write, for $I \subset \NZ$, ${\bf x}_I=  (x_t)_{t\in I} \in
[0,1]^{I}$.

Before stating precise assumptions on $p$, here is what one wants
to prove: $\exists C < \infty$, $m>0$ and a probability $\overline
p$ on $[0,1]$ such that $\forall E \subset [0,1]$, $E$ measurable,
\qq |p (x_T \in E | {\bf x}_{[-\infty,0]}) - {\overline p} (E)|
\leq C e^{-mT} \label{TM1} \qqq for all $T>0$ and all ${\bf
x}_{[-\infty,0]}$, where \qq p(x_T | {\bf x}_{[-\infty,0]}) =
\int^1_0 \prod^{T-1}_{t=1} dx_t \prod^{T}_{t=1} p({ x}_t | {\bf
x}_{[1,t-1]} \vee {\bf x}_{[-\infty,0]}) \label{TM2} \qqq and
${\bf x}_{[1,t-1]} \vee {\bf x}_{[-\infty,0]}$ denotes the obvious
configuration on $[-\infty, t-1]$.

Now let us state the assumptions on $p$ that will imply (\ref{TM1});
obviously, we assume
that:
\qq
p(x_t | {\bf x}_{[-\infty,t-1]}) \geq 0
\label{TM3}
\qqq
and
\qq
\int^1_0 dx_t p(x_t |{\bf x}_{[-\infty,t-1]}) = 1
\label{TM4}
\qqq
for all ${\bf x}_{[-\infty,t-1]}$.
Moreover, we assume that $p(\cdot | \cdot)$ is invariant under translations
of the
lattice $\NZ$, in a natural way. The non-trivial assumptions are:

\vs 2mm

a) Let, for $s<t-2$,
\qq
\delta_{s,t}({\bf x}_{[s,t]})  \equiv  p (x_t |
{\bf x}_{[s,t-1]} \vee {\bf 0}) - p(x_t| {\bf
x}_{[s+1,t-1]} \vee {\bf 0}),
\label{TM5}
\qqq
where ${\bf x}_I \vee {\bf 0}$ denotes the configuration equal to $x_t$ for
$t\in I$
and equal to zero elsewhere. We assume that  $\exists C < 0$, $m>0$ such that
$\forall s,t \in \NZ$ as above,
\qq
\| \delta_{s,t} \|_\infty \leq C \exp (-m|t-s|).
\label{TM6}
\qqq
\vs 2mm

b) Define, for $N\geq 1$, the Markov chain on $\Omega = [0,1]^N$ by
the transition probability
\qq
P ({\bf x}_{[1,N]}| {\bf x}_{[-N+1,0]}) =\prod^N_{t=1} p(x_t | {\bf
x}_{[t-N,t-1]}
\vee {\bf 0}).
\label{TM7}
\qqq

We assume that this Markov chain satisfies : $\exists \delta >0, \forall B
\subset
\Omega, \forall {\bf x}, {\bf x}' \in \Omega$,
\qq
P(B|{\bf x}) + P (B^c | {\bf x}') \geq \delta
\label{TM8}
\qqq
where $\delta$ is independent of $N$ (see however the Remark following the
proof of
Proposition 3.1 for a generalization).

\vs 3mm

\noindent
{\bf Proposition 3.1}. {\it Under assumptions a) and b) above, (\ref{TM1})
holds}.

\vs 3mm \no{\bf Remark 1.} The techniques used here can also prove
the analogue of (\ref{TM1}) with $x_T$ replaced by $x_{[T, T-L]}$,
for any finite $L$, and this, in turn, allows one  to associate to
the system of conditional probabilities a unique  probability
distribution on $[0, 1]^{\NZ}$ (which is called, in statistical
mechanics, the Gibbs state associated to the system of conditional
probabilities), but I will not go into that, because I want to
give here  only an elementary idea of the techniques used in
\cite{bk3}. Of course, this type of results is not new (see e.g.
\cite{Si}, Lecture 12, for a similar result, applied to dynamical
systems, with a somewhat different proof).

\vs 3mm

To prove the Proposition, we first use a result of Doob (\cite{Doob}, p.
197--198):

\vs 5mm

\no{\bf Lemma 3.1} {\it For the Markov chain defined in b) above, there
exists a
probability distribution $P$ on $\Omega$ such that $\forall {\bf x} \in
\Omega$,
$\forall B \subset \Omega$, $\forall n \geq 1$,}
\qq
|P^n (B | {\bf x}) - P (B) | \leq (1-\delta)^n.
\label{TM9}
\qqq

\vs 3mm

\no{\bf Proof.} Let ${\overline P} (B,n) = \sup_{{\bf x}} P^n (B |{\bf x})$ and
${\underline P}
(B,n) = \inf_{{\bf x}} P^n (B |{\bf x})$.
It is easy to see that ${\overline P} (B,n)$ is decreasing in $n$, while
${\underline P}
(B,n)$ is increasing in $n$. Thus, it
 is sufficient to prove the bound (\ref{TM9}) for the difference
$|\overline P
(B,n) - {\underline P} (B,n)|$ and, for that, we shall prove:
\qq
0 \leq {\overline P} (B, n+1) - {\underline P} (B, n+1) \leq (1-\delta)
({\overline P} (B,n)
- {\underline P}
(B,n)).
\label{TM10}
\qqq
Since $\overline P
(B,n) - {\underline P} (B,n)\leq 1$, (\ref{TM9}) follows.

Define a signed measure on subsets of $\Omega$:
\qq
\Psi_{{\bf x}, {\bf x'}} (E) = P (E|{\bf x}) - P (E|{\bf x}')
\label{TM11}
\qqq
and let $S^+$ (resp. $S^-$) denote the set where $\Psi_{{\bf x},{\bf x}'}
(E) \geq
0$ for $E \subset S^+$ (resp. $\leq 0$).

We have:
\qq
&&\overline P (B, n+1) - {\underline P} (B, n+1)
=\sup_{{\bf x},{\bf x}'} \int [P(d{\bf x}'' |
{\bf x}) - P(d{\bf x}'' | {\bf x}')] P^n (B|{\bf x''}) \non \\
&=& \sup_{{\bf x},{\bf x}'}
\int \Psi_{{\bf x},{\bf x}'} (d{\bf x}'') P^n (B|{\bf x}'') \non \\
&\leq& \sup_{{\bf
x},{\bf x}'}
(\Psi_{{\bf x},{\bf x}'} (S^+) \overline P (B,n) + \Psi_{{\bf x},{\bf x}'}
(S^-)
{\underline P}
(B,n)).
\label{TM12}
\en

By definition, $\Psi_{{\bf x},{\bf x}'} (S^-) = - \Psi_{{\bf x},{\bf x}'}
(S^+)$, so that
$$\Psi_{{\bf x},{\bf x}'} (S^+) \overline P (B,n) + \Psi_{{\bf x},{\bf x}'}
(S^-){\underline P}
(B,n)= \Psi_{{\bf x},{\bf x}'} (S^+) (\overline P (B,n) -{\underline P} (B,n)).
$$
Also, for any set $E \subset \Omega$, (\ref{TM8}) implies
$$
\Psi_{{\bf x},{\bf x}'} (E) = 1 -  (P(E^c |{\bf x}) + P(E|{\bf x}')) \leq
1-\delta.
$$
Applying this to $E=S^+$ in (\ref{TM12}) implies (\ref{TM10}).

\vs 3mm \no{\bf Remark 2.} We shall use this Lemma under the
following form: \qq \int_{\Omega} d{\bf x} |P^n({\bf x}|{\bf
x}')-P({\bf x})|\leq 2(1-\delta)^n \label{TM121} \qqq for all
${\bf x}'\in \Omega$; this follows by applying (\ref{TM9})
separately
 to the sets
where the integrand   is positive and negative.

Now, let us turn to the
\vs 3mm

\no{\bf Proof of Proposition 3.1.}

\vs 3mm
We write each factor in (\ref{TM2}) as
\qq
p(x_t| {\bf x}_{[-\infty,t-1]}) = p (x_t | {\bf x}_{[t-N,t-1]} \vee {\bf
0}) +
\sum_{|s-t|>N} \delta_{s,t}({\bf x}_{[s,t]}),
\label{TM13}
\qqq
where $N$ is an integer to  be chosen later. Insert this in the product in
(\ref{TM2}),
and expand: we get
\qq
p(x_T| {\bf x}_{[-\infty,0]}) = \sum_{I\subset [1,T]} \sum_{{\bf s}} \int
\prod^{T-1}_{t=1} dx_t \prod_{t\in I}\delta_{s,t}({\bf x}_{[s,t]})
\prod_{t\not\in
I} p(x_t | {\bf x}_{[t-N,t-1]} \vee {\bf 0}),
\label{TM14}
\qqq
where the sum over subsets $I$ corresponds to the choice in (\ref{TM13})
between
the first
term and the sum, while the sum over ${\bf s} = (s_t)_{t\in I}$ corresponds
to the
possible choices of
a term in that sum.

Now, let ${\overline I} =  \displaystyle{\bigcup_{t\in I}}   [s_t,t]$ and let
\qq
[1,T] \backslash {\overline I} = \bigcup_i J_i \bigcup_\alpha I_\alpha,
\label{TM15}
\qqq
where each $J_i$ is a union of intervals of length $N$, containing at least
two such
intervals, and each
$I_\alpha$ is an interval of length less than $2N$ between two connected
intervals in
$\overline I$ or  an interval of length less than $N$
between an interval in $\overline I$ and an interval $J_i$. The reason for
these definitions
is that, in the RHS of (\ref{TM14}), the only functions depending on $x_s$,
with $s$
in the complement of ${\overline I} $,
are the factors $p(x_t | {\bf x}_{[t-N,t-1]})$, so that, by integrating
over these variables, one can
obtain the transition probabilities of the Markov chain defined in
condition b) above.
For that, we need intervals of length at least $2N$, which are the $J_i$'s,
while the intervals
$I_\alpha$'s simply cover the leftover sites.

Since the model here is translation invariant, let us fix one interval $J_i
= J$, and
write it as a union of disjoint intervals of length
$N$: $J = \displaystyle{\bigcup^{n}_{l=0}} K_l$ with $K_l = [t +1 +lN,
t+(l+1)N]$.

We have, by definition (\ref{TM7}) of the transition probability $P$:
\qq
\int \prod^{n-1}_{l=1}  d{\bf x}_ {K_l}
\prod^{n}_{l=1} \prod_{t\in K_l} (p(x_t| {\bf
x}_{[t-N,t-1]}\vee
{\bf 0})) = P^n ({\bf x}_{K_{n}} | {\bf x}_{K_{0}}).
\label{TM16}
\qqq
Now write this as
\qq
P^{n} ({\bf x}_{K_{n}} | {\bf x}_{K_{0}}) - P({\bf x}_{K_{n}}) + P({\bf
x}_{K_{n}}),
\label{TM17}
\qqq
where $P$ is defined by (\ref{TM9}). Apply this to each interval $J_i$ in
(\ref{TM15}),
with
$n$ replaced by $n_i=\frac{|J_i|}{N}-1$. Insert that identity in (\ref{TM14})
for each
$J_i$ and expand the corresponding
product over $i$ of $A_i+B_i$, where $A=P^{n} ({\bf
x}_{K_{n}} | {\bf x}_{K_{0}}) - P({\bf x}_{K_{n}}) $ and $B= P({\bf
x}_{K_{n}})$.

For $E$ as in (\ref{TM1}), integrate over $E$ each term in the
resulting expansion, and  write
\qq
p(x_T \in E| {\bf x}_{[-\infty,0]}) = Q + R,
\label{TM18}
\qqq
where $Q$ collects all the terms in the resulting sum where at least one factor
$P({\bf x}_{K_{n_i}})$ appears and $R$ all the rest. Now,  the
presence of one such factor $P$ `decouples' $x_T$ from the initial
conditions ${\bf
x}_{[-\infty,0]}$, in the sense that, if we consider the difference
\qq
p(x_T\in E | {\bf x}_{[-\infty,0]}) - p (x_T \in E| {\bf x}'_{[-\infty,0]}),
\label{TM19}
\qqq
for two different past histories,
then the $Q$ sums are equal and only the $R$ sums contribute to the difference.
Indeed, fix a $K_{n_i}$ and consider all the terms in our expansions where
the factor
$P({\bf x}_{K_{n_i}})$ appears; let
$t_0$ be the last time before the interval
$K_{n_i}$. By construction, in all the terms under consideration, all the
functions that
depend on $x_t$, for $t>t_0$ do not depend on the variables $x_t$, for $t
\leq t_0$. So, if
we resum, in the expansion, all the terms depending on the
variables $x_t$, for $t \leq t_0$, we obtain, for the two terms in
(\ref{TM19}), $p(x_{t_0}|
{\bf x}_{[-\infty,0]})$, and
$p(x_{t_0}| {\bf x}'_{[-\infty,0]})$ (we simply use (\ref{TM14}) read from
right
to left, with $T$ replaced by
$t_0$). But performing in (\ref{TM2}) the integral over $x_{t_0}$ gives $1$
in both cases, which shows that
the difference between the respective
 sums cancel.

So, if we show that, $\exists C < \infty, m>0$ such that
\qq
|R| \leq C e^{-mT},
\label{TM20}
\qqq
$\forall {\bf x}_{[-\infty,0]}$, we obtain that the absolute value of
(\ref{TM19})
is exponentially small and, from that, (\ref{TM1}) easily follows.

Using the bound (\ref{TM6}) on $\delta_{s, t} $ and (\ref{TM121}) on
\qq
\int_{\Omega} \prod_{t\in K_{n}} dx_t | P^n ({\bf x}_{K_{n}} | {\bf
x}_{K_{0}}) - P
({\bf x}_{K_{n}})|,
\qqq
and the fact that, by (\ref{TM4}) and $x_t\in [0, 1]$,
all the integrals are bounded by $1$,
, we get:
\qq
|R| \leq \sum_I \sum_{{\bf s}} \prod_{t\in I} (Ce^{-m |t-s_t|}) \prod_i
(2(1-\delta)^{n_i}),
\label{TM21}
\qqq
where the second product runs over the intervals $J_i$ in
(\ref{TM15}), and
where $n_i = {|J_i|\over N} -1$. Note that the length of each $I_\alpha$ in
(\ref{TM15}) is less than $2N$ and, since such intervals are always
adjacent to a
connected component of $\overline I$ (unless $I=\emptyset$, in which case
this number is at
most $2$), the number of intervals
$I_\alpha$ is less than $2|I|+2$; the same bound holds for the
number of intervals $J_i$ in (\ref{TM15}) (in fact, a better bound holds here,
but we won't use it). So, we have:
\qq
\sum_i n_i\geq \sum_i \frac{|J_i|}{ N} -(2|I|+2) \geq \frac{(T-|{\overline
I}|)}{N}
-c|I|-2,
\label{TM211}
\qqq
for some number $c$, where, in the second inequality, we use $|I_\alpha|
\leq 2N$ and
(\ref{TM15}).

Using this, we can, by changing the constant $C$, bound (\ref{TM21}) by:
\qq
C\sum_I \sum_{{\bf s}} \prod_{t\in I} e^{-m|t-s_t|} C^{|I|}
(1-\delta)^{(T-|{\overline I}|)/N}
\label{TM22}
\qqq
Since by definition of ${\overline I}$, $\displaystyle{\sum_{t\in I}}
|t-s_t| \geq
|{\overline I}|$, we can, by considering separately
the terms where $|{\overline I}|\leq {T\over 2}$, and those where
 $|{\overline I}|> {T\over 2}$,  bound the sum in (\ref{TM22}) by
\qq
Ce^{- \tilde m T} \sum_I \sum_{{\bf s}} \prod_{t\in I} e^{-{m\over 2} |t-s_t|}
C^{|I|}
\label{TM23}
\qqq
where
\qq
\tilde m = \min \left({m\over 4}, {-\ln (1-\delta)\over 2N}\right).
\label{TM24}
\qqq

Now, choose $N$ so that
\qq
\sum_{|t-s| > N} e^{-{m\over 2} |t-s|} \leq \eta,
\label{TM25}
\qqq
with
\qq
(1+C\eta) \leq e^{\tilde m/2},
\label{TM26}
\qqq
which is possible since, from (\ref{TM25}) we see that, for large $N$,
$\eta =
\exp (-{\cal O}(N)) $ while, from ({\ref{TM24}), $\tilde m = {\cal O}
(N^{-1})$.

We use (\ref{TM25}) to control the sum over each $s_t$ in (\ref{TM23}), and
we get
\qq
(\ref{TM23}) \leq C e^{-\tilde m T} \sum_{I\subset [1,T]} (C \eta)^{|I|} \leq
C(1+C\eta)^T e^{-\tilde m T}
\label{TM27}
\qqq
and, using (\ref{TM26}), we get (\ref{TM20}) with $m = {\tilde m\over 2}$.

This completes the proof of  Proposition 3.1.

\vs 3mm \noindent{\bf Remark 3.} By considering (\ref{TM24},
\ref{TM25}, \ref{TM26}), we see that one can extend the proof to a
situation where $\delta$ in (\ref{TM8}) depends on $N$, as long
 as $\delta \geq \exp (-cN)$ for a constant $c$ small enough.

\vs 3mm
Now, let us turn to the real model, and make a list of the difficulties not
present
in the toy model. The first one is that time is continuous rather than
discrete,
but that is a minor problem. We can easily introduce a discretization of
time. A
more serious problem is that one deals with what are called ``unbounded
spins" in
statistical mechanics or what is also known as a ``large field problem",
namely the
variables $s(\tau)$ in (\ref{s1}), which play a role similar to the variables
$x_t$ here,
take value in $\NR^N$ rather than $[0,1]$ (actually, if we consider the
variable
$s$ over a unit time interval, they take values in a space of functions
from that
interval into $\NR^N$). And, what really causes a problem, is the fact that the
bounds (\ref{TM6}), (\ref{TM8}) do not hold when the variables $s$ take large
values. However, as we shall see in the next section, this is unprobable. Thus,
before doing an expansion as in (\ref{TM14}, \ref{TM17}), we must first
distinguish
between time intervals where the $s$ variables are large and those where
they are
small. Then, putting aside lots of technicalities, we perform the expansion
(\ref{TM14}) in the latter intervals and use estimates like (\ref{x1}) below to
control the sum over the intervals where $\omega$ is large.

Finally, there is an additional difficulty coming from the fact that the
definition
of the probabilities here involve a Girsanov representation. In statistical
mechanics, one usually deals with situations where the probabilities
(\ref{TM3})
can be written as:
\qq
p (x_t | {\bf x}_{[-\infty,t-1]}) = \exp (\sum_{t\in I} \phi_I (x_I)),
\label{TM28}
\qqq
where the  $\{\phi_I\}$'s represent ``many body interactions' (suitably
normalized so that (\ref{TM4}) holds) and
the sum runs over intervals $I\subset \NZ$ whose last point is $t$.
Then, a bound of the form
\qq
\| \phi_I \|_\infty \leq C\exp (-m  |I|),
\label{TM29}
\qqq
with $C<\infty$, $m>0$,
is enough to obtain (\ref{TM6}) and (\ref{TM8}). But here the probabilities
are not of that
form,
because of the stochastic integral $\int^t_0 f(\tau) \gamma^{-1} ds(\tau)$
in (\ref{s1}).

\section{A priori estimates on the transition probabilities}
\label{section 4}
\setzero\vskip-5mm \hspace{5mm }

The memory in the process (\ref{sred}) is coming from the
dependence of the solution of (\ref{leq}) on its initial
conditions. By Proposition 2.1, the dependence is weak if
$\int_0^t\|\Nna\omega\|^2$ is less than $cR t$ for a suitable $c$.
It is convenient to define, for each unit interval $ [n-1,n]\equiv
\Nn$, a quantity  measuring the size of $\omega$ on that interval
by: \qq D_n= \hf\sup_{t\in \Nn}
||\omega(t)||^{2}+\int_{\Nn}||\nabla\omega(t)||^{2}dt. \label{Dn}
\qqq

The following Proposition bounds the probability of the
unlikely event that we are interested in:

\vs 3mm

\noindent
{\bf Proposition 4.1}. {\it There
exist constants $c > 0$, $c'<\infty$,
 $\beta_0 < \infty$, such that
for all $t,t'$, $1 \leq t < t'$ and all $\beta\geq \beta_0$,
\qq
P \Bigl(\sum^{t'-1}_{n=t} D_{ n}
(\omega) \geq \beta R|t'-t| \Big|
\omega(0) \Bigr)
\leq  \exp({_1\over^R}c'e^{-t} \|\omega(0)\|^2)
\exp
(-c\beta |t'-t|).
\label{x1}
\qqq}
\vs 3mm

\vs 3mm

\noindent {\bf Remark 1.} This means that the probability that
$\omega$ is large over an interval of time decays exponentially
with the length of that interval, provided that $ \|\omega(0)\|$
is not too large. And, if $\|\omega(0)\|^2$ is of order $K$, $D_{
n} (\omega)$ will be, with large probability, of order $R$ after a
time of order $\log K$.

\vs 3mm

The main idea in the proof is
a probabilistic analogue of
the so-called  enstrophy balance: in the deterministic case,
using integration by parts
and $\Nna\cdot u=0$, on derives from (\ref{ve}) with $db=0$,
 the  identity:
$$
\frac{1}{2}\frac{d}{dt}\|\omega\|^2=-\|\nabla \omega\|^2,
$$
which implies that the enstrophy
($\| \omega\|^2$) decreases in time.
This basic property of  equation  (\ref{ve})
makes the proof of the following Lemma rather simple.

\vs 5mm

\no{\bf Lemma 4.1.} {\it For all $\omega(0)\in L^2$, and all $t\geq 0$,
\qq
E\Bigl[e^{{_{1}\over^{4R}}\|\omega(t)\|^2 }\;\Big|\;\omega(0)\Bigr]\leq
3e^{{_{1}\over^{4R}}e^{-t}\|\omega(0)\|^2},
\label{c4}
\qqq
and}
\qq
P( \|\omega(t)\|^2\geq D |\omega(0)) \leq 3
e^{-\frac{D}{4R}}
e^{{_{1}\over^{4R}}e^{-t}\|\omega(0)\|^2}.
\label{c5}
\qqq

\vs 3mm

\no{\bf Proof.} Let $x(\tau)= \lambda(\tau)\|\omega(\tau)\|^2 =
\lambda(\tau) \sum_k
|\omega_k|^2$ for $0\leq \tau\leq t $. Then by Ito's
formula (remember that, by (\ref{s3}), $\sum_k\gamma_k =
R$ and thus $\gamma_k\leq R$, $\forall k$):
\qq
{_d\over^{d\tau}}E[e^x]&=&
E[(\dot\lambda\lambda^{-1}x-2\lambda\sum_k k^2
|\omega_k|^2+{\lambda}
\sum_k\gamma_k+{2\lambda^2}\sum_k
\gamma_k|\omega_k|^2)e^x]\nonumber \\
&\leq& E[((\dot\lambda\lambda^{-1}-2+2\lambda
R)x+\lambda R)e^x],
\label{c6}
\qqq
where  $E$ denotes the conditional expectation, given $\omega(0)$,
and where we used the
Navier-Stokes
equation (\ref{ve}), $|k|\geq 1$ for $\omega_k\neq 0$, and the fact that the
nonlinear term does not contribute (using integration by parts
and $\Nna\cdot u=0$). Take now
$\lambda(\tau)= {_1\over^{4R}}e^{(\tau-t)}$ so that $\lambda \leq
{_1\over^{4R}}$, $\dot\lambda\lambda^{-1}= 1$,
$\dot\lambda\lambda^{-1}-2+2\lambda
R
\leq -\frac{1}{2}$ and $\lambda R\leq  \frac{1}{4}$. So,
\qq
{_d\over^{d\tau}}E[e^x]\leq E[(\frac{1}{4}-\frac{1}{2}x)e^x]\leq
\frac{1}{2}-\frac{1}{4}E[e^x],
\non
\qqq
where the last inequality follows by using
$(1-2x)e^x\leq
2-e^x$. Thus,  Gronwall's
inequality implies that:
\qq
E[e^{x(\tau)}]\leq e^{-\frac{\tau}{4}}e^{x(0)} +2\leq 3e^{x(0)},
\non
\qqq
i.e., using the definition of $\lambda(\tau)$,
\qq
E\Bigl[\exp(\frac{e^{\tau-t}}{4R}\|\omega(\tau)\|^2 )\Bigr]\leq
3\exp(\frac{e^{-t}\|\omega(0)\|^2}{4R}).
\non
\qqq
This proves (\ref{c4}) by putting $\tau=t$;
(\ref{c5}) follows from (\ref{c4}) by Chebychev's
inequality.

\vs 3mm

Since the $D_n$ in (\ref{x1}) is
the supremum over unit time intervals of
\qq
D_t(\omega) = {1\over 2} \| \omega (t) \|^2 + \int^t_{n-1}
\|
\Nna \omega \|^2 d \tau \hspace{5mm} n-1 \leq t \leq n,
\label{a0}
\qqq
which does not involve only
 $ \| \omega (t) \|^2$, we need
to control also the evolution of $D_t(\omega)$ over a
unit time interval, taken, for now,
to be $[0,1]$. From the Navier-Stokes equation
(\ref{ve}) and Ito's formula, we obtain
\qq
D_t(\omega) = D_0 (\omega) + Rt + \int^t_0 (\omega, db)
\label{a01}
\qqq
(since the nonlinear term does not contribute, as in (\ref{c6})).

Our basic estimate is:

\vs 3mm

\no{\bf Lemma 4.2.} {\it There exist $C < \infty$, $c>0$ such that,
$\forall A\geq 3D_0(\omega)$
\qq
P(\sup_{t\in [0,1]} D_t (\omega) \geq  A|\omega(0) )
\leq C \exp(-{_{cA}\over^{R}}).
\label{a1}
\qqq
}
\vs 3mm

\noindent {\bf Remark 2.} While the previous Lemma showed that
$\|\omega (t) \|^2$ tends to decrease as long as it is larger than
${\cal O} (R)$, this Lemma shows that, in a unit interval, $D_t
(\omega)$ does not increase too much relative to $D_0 (\omega) =
{1\over 2} \|\omega (0)\|^2$. Thus, by combining these two Lemmas,
we see that $D_n (\omega)=\displaystyle{\sup_{t\in[n-1,n]}} D_t
(\omega)$ is, with large probability, less than $\|\omega
(0)\|^2$, when the latter is larger than ${\cal O} (R)$, at least
for $n\geq n_0$ not too small. Thus, it is unlikely that $D_n
(\omega)$ remains much larger than $R$ over some interval of
(integer) times, and this is the basis of the proof of Proposition
4.1.

\vs 3mm

Without entering
into details, here are the main ideas in  the proof of (\ref{a1}).
From (\ref{a01}), we see that it is enough to get an upper bound on
\qq
 P\left(\sup_{t\in [0,1]} |
\int^t_0 (\omega,
db)| \geq (A - D_0 - R)\Big|\omega(0)\right).
\label{PE}
\qqq
We use (see \cite{bk3} for more details) Doob's inequality (see
e.g.\cite{simon}, p.24), to
reduce the control over the supremum over $t$
to  estimates on $|\int^1_0 (\omega, db)|$.
Letting $E$ denote the conditional expectation, given $\omega(0)$.
and using Novikov's bound (see e.g. the proof of Lemma 5.2 below), we get
\qq
&& E(e^{\pm \varepsilon \int^1_0 (\omega, db)})
 \leq \left( E
(e^{2\varepsilon^2\int^1_0 d\tau(\omega(\tau),\gamma
\omega(\tau))})\right)^{1/2}\non\\
&\leq& \left(\int^1_0 d\tau E (e^{2\varepsilon^2
(\omega(\tau),\gamma \omega(\tau))})\right)^{1/2} \leq
\left(\int^1_0 d\tau E (e^{2\varepsilon^2 R \| \omega
(\tau)\|^2})\right)^{1/2} \label{a04} \qqq where the last two
inequalities follow from Jensen's inequality, applied to
\linebreak $e^{2\varepsilon^2\int^1_0 d\tau(\omega(\tau),\gamma
\omega(\tau))}$, and  from $\gamma_k \leq  R$ (see (\ref{s3})).
Now, choosing $\varepsilon$ so that
$2\varepsilon^2R=\frac{1}{4R}$, i.e. $\varepsilon =\frac{1}{{\sqrt
8} R}$, we can use (\ref{c5}) to bound the RHS of (\ref{a04}).
Combining this with Chebychev's inequality gives bounds on
(\ref{PE}).

\section{Decoupling estimates}
\label{section 5}
\setzero\vskip-5mm \hspace{5mm }

In this section, I shall give a very
brief sketch of the ideas used to prove the analogue
of assumptions (\ref{TM6}) and (\ref{TM8}) of section 3 in the present
setting, at least
in the probable regions where $\omega$ is small. The main point is to
understand
the analogue of the bound (\ref{TM29}), which expresses the exponential
decay of interactions.
What plays the role of the right hand side of (\ref{TM28}) is, see (\ref{s1}):
\qq
g_t\equiv e^{\int_{t-1}^t(f(\tau),\gamma^{-1}(ds(\tau)
-\hf f(\tau)d\tau))}
\qqq
where, for simplicity, I consider a unit time interval $[t-1, t]$.
We want to show that this depends weakly on the past; so consider two
functions $g_1$, $g_2$,
defined in terms of two functions $f_1$, $f_2$, themselves defined through
different
$l_1$ and $l_2$
(see (\ref{f}, \ref{oomega})). And, by analogy with what we did in section
3, we choose
$l_1= l(t,{\bf s}([0,t]),l(0)=0)$, $l_2= l(t,{\bf s}([1,t]),l(1)=0)$, i.e.
we set
the large $k$ modes equal to zero at different times ($0$ or $1$).
Using (\ref{22}), we see that
$l_1= l(t,{\bf s}([1,t]),l_1(1))$, with $l_1(1)= l_1(1,{\bf
s}([0,1]),l(0)=0)$, so that we
have, at time $t=1$, two initial conditions, $l_1(1)$, $l_2(1)=0$, with
$\|l_1(1)-l_2(1)\|=\|l_1(1)\|$ of order one, if $\omega$ is small in
 the interval $[0,1]$.

Now, if $g_t$
depends weakly on the past, it should mean that, for large $t$,
 $g_1$ and $g_2$ are, in some sense,
exponentially close. To measure the difference, write:
\qq
g_1-g_2=(1-\frac{g_2}{g_1})g_1\equiv (1-H)g_1,
\qqq
which will be convenient, since we deal
 with unbounded variables for which sup norm estimates like in
(\ref{TM29}) are not available. Explicitly: \qq
H=e^{\int_{t-1}^t(\delta f(t),\gamma^{-1}(ds(t)- f_{1}(t)dt)) -\hf
\int_{t-1}^t(\delta f(t),\gamma^{-1}\delta f(t))dt } \label{DE1}
\qqq where $\delta f= f_2-f_1$. What we want to show is that $1-H$
is, in a suitable sense, small.

The next Lemma gives a bound on $\|\delta f\|$ in terms of $\|\delta l\|$,
and $\| \omega
\|$;
$\|\delta l\|$ is controlled by Proposition 2.1, provided that $\omega$
is small, in the sense discussed in section 4, in which case
$\| \omega
\|$ is also controlled, using $\sup _{t\in \Nn}\| \omega(t)
\|
\leq (2D_n)^{\hf}$.
\vs 3mm

\noindent {\bf Lemma 5.1}.
{\it Let $f(\omega)=PF(\omega)$ and $\omega=s+l$,
$\omega'=s+l'$. Then,
\qq
\| \delta f\| = \|f (\omega) - f (\omega')\| \leq
C(R)(2\| \omega \| \| \delta l \| + \|
\delta l \|^2)
\label{a13}
\qqq
with $\delta l = l - l'$ and $C(R)$ a constant depending
 on the parameter $R$ (see (\ref{s3})).}

\vs 3mm

\noindent
{\bf Proof}. We have
\qq
|f_k (\omega) - f_k (\omega')|
\leq \sum_p |\omega_{k-p} \omega_{p} -
\omega'_{\kappa-p} \omega'_{p} | {|k|\over |p|}
\non
\qqq
which, since $|k| \leq \sqrt{\kappa R}$ is bounded by
\qq
\sqrt{\kappa R} \sum_p | s_{k-p}
\delta l_p + s_p \delta l_{k-p} + l_p
l_{k-p} - l'_p l'_{k-p}|.
\label{a15}
\qqq
Writing $l_p l_{k-p} - l'_p l'_{k-p} =
l_p \delta l_{k-p} + l_{k-p} \delta
l_p - \delta l_p \delta l_{k-p}$ and using Schwarz'
inequality,  we get
\qq
(\ref{a15}) \leq \sqrt{\kappa R} (2 \| \omega \| \| \delta l
\| \| + \|\delta l \|^2)
\non
\qqq
which proves (\ref{a13}), since $f_k \neq 0$ only for
$k\leq \kappa R$, so that the sum in the $L^2$ norm $\| \delta f\|$
 runs
over $C(R)$ terms.

This Lemma would be enough to control $(1-H)$  if we had only
in (\ref{DE1}) the factor $e^{
-\hf \int_{t-1}^t(\delta f(t),\gamma^{-1}\delta f(t))dt
}$, which involves only ordinary integrals.

To control the stochastic integral, it is convenient
to  undo the Girsanov transformation, i.e. to change
variables from $s$ back to $b$. Let $E$ denote the
expectation with respect to the Brownian motion $b$ with
covariance $\gamma$ on the time interval $[t-1, t]$. We get, using
(\ref{sred}):
\qq
H=e^{\int_{t-1}^t(\delta f(t),\gamma^{-1}db(t))
-\hf\int_{t-1}^t(\delta f(t),\gamma^{-1}\delta f(t))dt
}.
\label{dg11}
\qqq

Write now $(1-{H})^2=1-2{ H} +{ H}^2$; to give a flavour of the estimates,
let us see
how one could show that the expectation with respect to $E$
of $-2{ H} +{ H}^2$is close to $1$, i.e. that the expectation
of $(1-H)^2$ is close to zero.
One can rather easily bound from below
the expectation of $H$, using Jensen's inequality;
to get an upper bound on the expectation of ${ H}^2$,
one uses:

\vspace*{2mm}

\no{\bf Lemma 5.2.} {\it Let $\zeta(t)\in C([0,1], H_s)$ be
progressively
measurable. Then
\qq
Ee^{\int_0^1 (\zeta,\gamma^{-1}db)
+\lambda\int_0^1(\zeta,\gamma^{-1}\zeta) dt}\leq e^{2(1+\lambda)
||\zeta||^2\rho^{-1}}
\label{novikov}
\qqq
where $||\zeta||=\sup_\tau||\zeta(\tau)||_2$.}

\vs 2mm

\no{\bf Proof}. This is just a Novikov bound: we bound the LHS,
using Schwarz' inequality, by
$$
(Ee^{\int_0^1 (2\zeta,\gamma^{-1}db)
-2\int_0^1(\zeta,\gamma^{-1}\zeta dt)})^\hf
(Ee^{2(1+\lambda)\int_0^1(\zeta,\gamma^{-1}\zeta) dt})^\hf
$$
and note that the expression inside the first square root
is the expectation of a martingale and equals one.

\vs 2mm

We can then apply this Lemma to $\zeta =2\delta f$, $\lambda=-\frac{1}{4}$,
replacing $[0, 1]$ by $[t-1, t]$,
and  use the estimates
coming from Lemma 5.1 and Proposition 2.1 to show that the RHS of
(\ref{novikov})
is exponentially close to $1$, for $t$ large. This gives a rough idea of
why the ``interactions" here are
exponentially decaying, but it must be said that the full story
is far more complicated and I refer to reader
to \cite{bk3} for more details.

\vspace*{2mm}

\no{\bf Acknowledgments.}

I wish to thank my coauthors, Antti Kupiainen and Raphael Lefevere
without whom this work would not have been possible.
This work was supported in part by ESF/PRODYN.

\end{document}